\documentclass[10pt]{amsart}
\usepackage{amsmath,amsthm, amscd, amsfonts, amssymb, graphicx, xcolor}
\usepackage[colorlinks,linkcolor=blue,urlcolor=blue,citecolor=blue,hyperindex,dvipdfm]{hyperref}
\newtheorem{thm}{Theorem}[section]
\newtheorem{cor}[thm]{Corollary}
\newtheorem{lem}[thm]{Lemma}

\newtheorem{exam}[thm]{Example}
\numberwithin{equation}{section}

\begin{document}

\title{The Generalized Flanders' Theorem in Unit-regular Rings}

\author{Dayong Liu}
\author{Aixiang Fang$^*$}
\address{
College of Science \\ Central South University of Forestry and Technology  \\ Changsha , China}
\email{<liudy@csuft.edu.cn>}
\address{
College of Mathematics and Physics \\ Hunan University of Arts and Science \\ Changde,  China}
\email{<fangaixiangwenli@163.com>}

\thanks{$^*$Corresponding author}

\begin{abstract}
Let $R$ be a unit-regular ring, and let $a,b,c\in R$ satisfy $aba=aca$. If $ac$ and $ba$ are group invertible,
we prove that $ac$ is similar to $ba$. Furthermore, if $ac$ and $ba$ are Drazin invertible,
then their Drazin inverses are similar. For any $n\times n$ complex matrices $A,B,C$ with $ABA=ACA$ , we prove that $AC$ and $BA$ are similar
if and only if their $k$-powers have the same rank. These generalize the known Flanders' theorem proved by Hartwig.
 \end{abstract}

\maketitle

\section{Introduction}

An element $a\in R$ has Drazin inverse if there exists an element $x\in R$ such that
\begin{equation}\nonumber
 ax=xa, \ \ xax=x, \ \ a^{k+1}x=a^{k} \   for \ some \  k\in {\mathbb{N}},
\end{equation}
 or equivalently,
\begin{equation}\nonumber
 ax=xa, \ \ xax=x, \ \ a^2x-a \in N(R)
\end{equation}
where $N(R)$ denotes the set of all nilpotents in $R$. If $a$ is Drazin invertible, the Drazin inverse of $a$ is unique,  denote $x$ by $a^D$.
The least nonnegative $k$ which satisfies formulas above is called the index of $a$, denoted by $ind(a)$.
If $ind(a)$ =1, $a$ is said to be group invertible. In this case,  the element $x$ is called the group inverse of $a$ and denoted by $a^{\#}$, that is ,
\begin{equation}\nonumber
 aa^{\#}=a^{\#}a, \ \ a^{\#}aa^{\#}=a^{\#}, \ \ aa^{\#}a=a.
\end{equation}
\indent
We use $R^{\#}$ to stand for the set of all group invertible elements of $R$.   \\
\indent
Two elements $a, b \in R$  are similar, i.e., $a\sim b$, if there exists an invertible element $s$ such that $a=s^{-1}bs$.  \\
\indent
In~\cite{HARTWIG}, Hartwig has proved the Flanders' theorem, and given a generalization, i.e., a strongly $\pi$-regular ring $R$ is unit-regular if and only if $R$ is regular and the Drazin inverses of $ab$ and $ba$ are similar.  \\
\indent
For any $A,B\in R^{n\times n}$ over a Bezout domain $R$, Cao and Li~\cite[Theorem 3.6]{CAO} proved that
\begin{equation}\nonumber
 (AB)^{\#} \ and \ (BA)^{\#} \ exist, \ then \ AB \sim BA.
\end{equation}
In 2016, Deng~\cite[Theorem 2.6]{DENG} generalized this results to the operators on an arbitrary Hilbert space.
Furthermore, Mihallovi$\mathrm{\acute{c}}$  and Djordjevi$\mathrm{\acute{c}}$~\cite[Theorem 2.2]{MIHAJLOVIC} claimed the preceding results for $a,b$ in a general ring.
Unfortunately, there is a gap in its proof. To be specific, for any $a,b\in R$, generally, the conditions that $(ab)^{\#}$ and $(ba)^{\#}$ exist do not imply that $ab \sim ba$.
We will give a counter-example in the next section.  \\
\indent
Recall that a ring $R$ is unit-regular provided that for each $a\in R$, there is a unit $u\in R$ such that $aua=a$.
For example, the complex field ${\mathbb{C}}$ is unit-regular.
In this paper, we proved that for a unit-regular ring $R$, $a,b,c\in R$ with $aba=aca$, if $(ac)^{\#}$ and $(ba)^{\#}$ exist, then $ac\sim ba$ and $(ac)^{\#}\sim (ba)^{\#}$.
We further generalized the Flanders' theorem and proved that for $A,B,C\in {\mathbb{C}}^{n\times n}$ with $ABA = ACA$, then $AC\sim BA$ if and only if $rank(AC)^k = rank(BA)^k$ for $k=1,2,\cdots$. Moreover, if $s\geq max\{ind(AC), ind(BA)\}$, then $(AC)^s\sim (BA)^s$.    \\
\indent
Throughout this paper, all rings are associative with an identity, the set of all invertible elements of $R$ will be denoted by $U(R)$. ${\mathbb{N}}$ stands for the set of all natural numbers. \\

\vskip1cm

\section{Main Results}

We begin with a counter-example which infers that~\cite[Theorem 2.2]{MIHAJLOVIC} is not true.
\begin{exam}\rm
Let $V$ be an infinite dimensional vector space of a field $\mathbb{F}$, and let $R=End_{\mathbb{F}}(V)$.
Let $\{x_1, x_2, \cdots, x_n, \cdots\}$ be a basis of $V$.\\
\indent
Definition
\begin{center}
 $\sigma(x_i)=x_{i+1}$ for all $i\in {\mathbb{N}}$,   \\
 $\tau(x_1)=0, \ \ \tau(x_i)=x_{i-1}$ for all $ i\geq 2 $.
\end{center}
Then $\sigma, \tau \in R$, and for any $i\in {\mathbb{N}}$,
$$ \tau\sigma(x_i)=\tau(x_{i+1})=x_i, $$
i.e., $\tau\sigma = 1_{V}$. Therefore $\tau\sigma$ is invertible in $R$, hence $\tau\sigma\in R^{\#}$.  \\
\indent
Since $\sigma\tau(x_1)=\sigma(0)=0$, we have $\sigma\tau \neq 1$. But $(\sigma\tau)^2=\sigma(\tau\sigma)\tau=\sigma\tau$, $\sigma\tau$ is an idempotent, and so $\sigma\tau\in R^{\#}$. \\
\indent
We claim that $\tau\sigma \nsim \sigma\tau$, otherwise, there exists $s\in U(R)$ such that
$\tau\sigma\cdot s=s \cdot \sigma\tau$ which implies $\sigma\tau=1$, a contradiction.  \hfill$\Box$
\end{exam}
\begin{lem}\label{lem1}Let $a,b,c\in R$ satisfy $aba=aca$.     \\
\indent
$\mathrm{(1)}$   If  $(ac)^{D}$ or $(ba)^{D}$ exists, then
$$(ba)^{D}=b[(ac)^{D}]^2 a,  \ \ (ac)^{D}=a[(ba)^{D}]^2 c,$$
$$ a(ba)^{D}=(ac)^{D}a, \ \ ab(ac)^D=ac(ac)^D.$$
\indent
$\mathrm{(2)}$  If $(ac)^{\#}$ and $(ba)^{\#}$ exist, then
$$(ba)^{\#}=b[(ac)^{\#}]^2 a, \ \ (ac)^{\#}=a[(ba)^{\#}]^2 c, $$
$$ab(ac)^{\#}=ac(ac)^{\#}, \ \ a(ba)^{\#}=(ac)^{\#}a. $$
\end{lem}
\noindent
{\it Proof}. (1) In view of~\cite[Theorem 2.7]{ZENG}, we have $(ba)^{D}=b[(ac)^{D}]^2 a$ and $(ac)^{D}=a[(ba)^{D}]^2 c$.
Moreover, we assume that $(ac)^D$ exists, we get
$$a(ba)^D =ab[(ac)^D]^2 a=ac[(ac)^D]^2 a=(ac)^D a,$$
$$ab(ac)^D = abac[(ac)^D]^2=acac[(ac)^D]^2 =ac(ac)^D. $$
\indent
(2) Suppose that $(ac)^{\#}$ and $(ba)^{\#}$ exist. Then $(ac)^{\#} = (ac)^D$ and $(ba)^{\#}=(ba)^D$, we obtain the result by (1).   \hfill$\Box$        \\
\indent
We are ready to prove:  \\
\begin{thm}\label{thm1}
Let $R$ be a unit-regular ring and let $a, b, c\in R$ satisfy $aba=aca$. If $(ac)^{\#}$ and $(ba)^{\#}$ exist, then $ac \sim ba$.
\end{thm}
\noindent
{\it Proof}. Choose $x=b(ac)^{\#}$, $y=ac(ac)^{\#}a$. Since $aba=aca$, we check that
$$\begin{array}{rcl}
x(ac)y  & = &  b(ac)^{\#} \cdot ac \cdot ac(ac)^{\#}a=babac[(ac)^{\#}]^2 a=bab[(ac)^{\#}]^2 aba=ba;  \\
y(ba)x  & = &  ac(ac)^{\#}a \cdot ba \cdot b(ac)^{\#}=ac(ac)^{\#}acac(ac)^{\#}=ac;  \\
xyx  & = &  b(ac)^{\#} \cdot ac(ac)^{\#}a \cdot b(ac)^{\#}=b(ac)^{\#}ac(ac)^{\#}=x;  \\
yxy  & = &  ac(ac)^{\#}a \cdot b(ac)^{\#} \cdot ac(ac)^{\#}a=(ac)^{\#}acac(ac)^{\#}a=y.
\end{array}$$
\noindent
Since $R$ is unit-regular, we have $x=xvx$ for some $v\in U(R)$.
Set $$u=(1-xy-xv)v^{-1}(1-yx-vx).$$
Since $(1-yx-vx)^2=1$ and $(1-xy-xv)^2=1$,
we get
$$(1-xy-xv)v^{-1}(1-yx-vx)^2 v(1-xy-xv) =1,$$
$$(1-yx-vx)v(1-xy-xv)^2 v^{-1}(1-yx-vx)=1,$$
i.e., $u$ is invertible in $R$. Furthermore, we have
$$u^{-1}=(1-yx-vx)v(1-xy-xv)=v-vxv+y.$$
Then we verify that
$$(ac)u^{-1}=ac(v-vxv+y)=y(ba)xv(1-xv)+acac(ac)^{\#}a=aca,$$
$$u^{-1}(ba)=(v-vxv+y)ba=(1-vx)vx(ac)y+ac(ac)^{\#}aba=aca.$$
Therefore $(ac)u^{-1}= u^{-1}(ba)$. Hence $ac=u^{-1}(ba)u$, as desired.       \hfill$\Box$
\begin{cor}\label{cor1}
Let $ A , B , C\in \mathbb{C}^{n\times n}$ with $ABA=ACA$. If $(AC)^{\#}$ and $(BA)^{\#}$ exist, then $ AC \sim BA $.
\end{cor}
\noindent
{\it Proof}. By~\cite[Corollary 4.5]{CHEN2010}, the matrix ring ${\mathbb{C}^{n\times n}}$ is unit-regular.
The result follows by the Theorem~\ref{thm1}.       \hfill$\Box$      \\
\vskip-2mm
\indent
We turn to consider the similarity of the group inverses.
\begin{thm}\label{thm2}
Let $R$ be a unit-regular ring and let $a, b, c\in R$ with $aba=aca$. If $(ac)^{\#}$ and $(ba)^{\#}$ exist,
then $(ac)^{\#} \sim (ba)^{\#}$.
\end{thm}
\noindent
{\it Proof}. Similarly to Theorem~\ref{thm1}, let $x=b(ac)^{\#}$, $y=ac(ac)^{\#}a$. We can verify the following formulas:
$$x(ac)^{\#}y=(ba)^{\#}, \ \ \ y(ba)^{\#}x=(ac)^{\#}, \ \ \ xyx=x, \ \ \ yxy=y.$$
\noindent
Since $R$ is unit-regular, we have $x=xvx$ for some $v\in U(R)$. Set
$$ u=(1-xy-xv)v^{-1}(1-yx-vx). $$
\noindent
Analogously to Theorem~\ref{thm1}, we have
$$u^{-1}=(1-yx-vx)v(1-xy-xv)=v-vxv+y. $$
Furthermore, as in the proof of Theorem~\ref{thm1}, we can check that
$$(ac)^{\#}u^{-1}=y(ba)^{\#}xv(1-xv)+(ac)^{\#}ac(ac)^{\#}a=(ac)^{\#}a,$$
$$u^{-1}(ba)^{\#}=(1-vx)vx(ac)^{\#}y+ac(ac)^{\#}a(ba)^{\#}=(ac)^{\#}a.$$
Thus we get $(ac)^{\#}=u^{-1}(ba)^{\#}u $. This completes the proof.       \hfill$\Box$
\begin{cor}\label{cor2}
Let $A, B, C\in \mathbb{C}^{n\times n}$ with $ABA=ACA$. If $(AC)^{\#}$ and $(BA)^{\#}$ exist,
then $(AC)^{\#} \sim (BA)^{\#}$.
\end{cor}
\noindent
{\it Proof}. Analogously to Corollary~\ref{cor1}, $\mathbb{C}^{n\times n}$ is unit-regular.
Therefore we complete the proof by Theorem~\ref{thm2}.       \hfill$\Box$
\begin{exam} \rm
Let
 $A=\begin{bmatrix}
     0 & 1 \\
     0 & 1 \\
   \end{bmatrix},
B=\begin{bmatrix}
     1 & 1 \\
     1 & 1 \\
   \end{bmatrix},
C=\begin{bmatrix}
     0 & 0 \\
     1 & 1 \\
   \end{bmatrix}
\in \mathbb{C}^{2\times 2}$. Then $B\neq C$ and $ABA=ACA$.
$AC=\begin{bmatrix}
     1 & 1 \\
     1 & 1 \\
   \end{bmatrix}$,
$BA=\begin{bmatrix}
     0 & 2 \\
     0 & 2 \\
   \end{bmatrix}$.
It is easy to verify that
$$(AC)^{\#}=
\begin{bmatrix}
    \frac{1}{4} & \frac{1}{4} \\
     \frac{1}{4} & \frac{1}{4} \\
   \end{bmatrix},
(BA)^{\#}=
\begin{bmatrix}
    0 & \frac{1}{2} \\
    0 & \frac{1}{2} \\
   \end{bmatrix}.$$
In view of Corollary~\ref{cor1} and~\ref{cor2}, $AC\sim BA$ and $(AC)^{\#}\sim (BA)^{\#}$.
Actually, let
$S=\begin{bmatrix}
     1          &      0      \\
    \frac{1}{2} & \frac{1}{2} \\
   \end{bmatrix}$,
then $AC=S^{-1}BAS$ and $(AC)^{\#}=S^{-1}(BA)^{\#}S$.
Therefore $AC\sim BA$ and $(AC)^{\#} \sim (BA)^{\#}$.        \hfill$\Box$
\end{exam}
\indent
Hartwig~\cite[Theorem 1]{HARTWIG} proved that $x^D \sim y^D$ if and only if $x^2 x^D \sim y^2 y^D$ over a strongly $\pi$-regular ring $R$ .
 We will give a generalization of this result as following.
\begin{thm}\label{thm3}
Let $R$ be a unit-regular ring and let $a, b, c\in R$ with $aba=aca$. If $(ac)^{D}$ or $(ba)^{D}$ exists, then $(ac)^{D}\sim(ba)^{D}$.
In this case,  $(ac)^2(ac)^D\sim(ba)^2(ba)^D$.
\end{thm}
\noindent
{\it Proof}. By virtue of Lemma~\ref{lem1}, we have $(ba)^{D}=b[(ac)^D]^2 a = b(ac)^D \cdot (ac)^D a $. Let
$$x=b(ac)^D, \ \ y=ac(ac)^D a. $$
Then we check that
$$\begin{array}{rcl}
x(ac)^D y  &  =  &  b(ac)^D \cdot (ac)^D \cdot ac(ac)^D a =(ba)^D;                             \\
y(ba)^D x  &  =  &  ac(ac)^D a \cdot (ba)^D \cdot b(ac)^D=ac(ac)^D (ac)^D ac(ac)^D=(ac)^D;     \\
      xyx  &  =  &  b(ac)^D \cdot ac(ac)^D a \cdot b(ac)^D=b(ac)^D ac (ac)^D=x;                \\
      yxy  &  =  &  ac(ac)^D a \cdot b(ac)^D \cdot ac(ac)^D a=ac(ac)^D ac(ac)^D a=y;
\end{array}$$
By the unit-regularity of $R$, we have $x=xvx$ for some $v\in U(R)$. Set
$$ u=(1-xy-xv)v^{-1}(1-yx-vx).$$
\noindent
Then
$$ u^{-1}=(1-yx-vx)v(1-xy-xv)=v-vxv+y. $$
We verify that
$$(ac)^D u^{-1}=y(ba)^D xv(1-xv)+(ac)^D ac(ac)^D a=(ac)^D a,$$
$$u^{-1}(ba)^D=(1-vx)vx(ac)^D y+ac(ac)^D a(ba)^D=(ac)^D a.$$
Therefore
$$(ac)^D=u^{-1}(ba)^D u.$$
i.e., $(ac)^{D}\sim (ba)^{D}$.  \\
\indent
Accordingly, by~\cite[Theorem 1]{HARTWIG}, $(ac)^2 (ac)^D \sim (ba)^2 (ba)^D $.        \hfill$\Box$
\begin{cor}\label{cor3}
Let $A, B, C\in\mathbb{C}^{n\times n}$ with $ABA=ACA$. Then $AC\sim BA$ if and only if
$rank(AC)^k = rank(BA)^k$ for $k=1,2,\cdots$.
\end{cor}
\noindent
{\it Proof}. The necessity is obvious. For the sufficiency, there exist two invertible matrices $P$ and $Q$ such that
\begin{equation}\nonumber
PACP^{-1}=\left(
\begin{array}{cc}
U_{1}  &    \\
       &  N_{1}
\end{array}\right) , \ \
QBAQ^{-1}=\left(
\begin{array}{cc}
U_{2}  &    \\
       &  N_{2}
\end{array}\right)
\end{equation}
where $U_1, U_2$ are invertible, and $N_1, N_2$ are nilpotent.  Choose $s=ind(N_1)+ind(N_2)$, then
\begin{equation}\nonumber
P(AC)^s P^{-1}=\left(
\begin{array}{cc}
U_{1}^s  &    \\
         &  O
\end{array}\right), \ \
Q(BA)^s Q^{-1}=\left(
\begin{array}{cc}
U_{2}^s  &    \\
         &   O
\end{array}\right).
\end{equation}
Since $rank(AC)^s = rank(BA)^s$, we see that $U_1$ and $U_2$ have the same rank.
It is easy to check that
\begin{equation}\nonumber
(AC)^D=P^{-1}\left(
\begin{array}{cc}
U_{1}^{-1}  &    \\
            &  O
\end{array}\right)P, \ \
(BA)^D=Q^{-1}\left(
\begin{array}{cc}
U_{2}^{-1}  &    \\
            &  O
\end{array}\right)Q.
\end{equation}
In view of Theorem~\ref{thm3}, we have $$(AC)^D \sim (BA)^D,$$
i.e.,
\begin{equation}\nonumber
\left(
\begin{array}{cc}
U_1^{-1}  &        \\
          &   O
\end{array}
\right)
\sim
\left(
\begin{array}{cc}
U_2^{-1}  &        \\
          &   O
\end{array}
\right),
\end{equation}
which follows $U_1^{-1}\sim U_2^{-1}$, and so $U_1\sim U_2$.
Moreover, as $rank(AC)^k = rank(BA)^k$,  $rank(N_1)^k = rank(N_2)^k$ for all positive integers $k$.
Since $N_1, N_2$ are nilpotent, by the Jordan forms of $N_1, N_2$, we have that $N_1$ and $N_2$ have the same Jordan forms. Hence $N_1 \sim N_2$. Therefore $AC\sim BA$, as asserted .       \hfill$\Box$
\begin{thm}\label{thm4}
Let $A, B, C\in {\mathbb{C}}^{n\times n}$ with $ABA=ACA$. If $s\geq max\{ind(AC)$, $ind(BA)\}$, then $(AC)^s \sim (BA)^s $.
\end{thm}
\noindent
{\it Proof}. For $s\geq max\{ind(AC), ind(BA)\}$, we shall show that $(AC)^s$ and $(BA)^s$ have group inverses. Obviously, we have
$$(AC)^s [(AC)^D]^s  = AC(AC)^D , \ \ \  [(AC)^D]^s (AC)^s  = AC(AC)^D, $$
i.e.,
$$(AC)^s[(AC)^D]^s=[(AC)^D]^s (AC)^s.$$
Since $AC(AC)^D$ and $(AC)^D AC$ are idempotent, one can check that
$$[(AC)^D]^s (AC)^s [(AC)^D]^s = AC(AC)^D[(AC)^D]^s = [(AC)^D]^s,  $$
$$(AC)^s [(AC)^D]^s (AC)^s = AC (AC)^D(AC)^s = (AC)^s.  $$
So $(AC)^s$ is group invertible.  Similarly, $(BA)^s$ is also group invertible.  \\
\indent
Let $C'=C(AC)^{s-1}$, $B'=(BA)^{s-1}B$, then $AC'A=AB'A$.
Then apply Theorem~\ref{thm1} to our case, we obtain the result.       \hfill$\Box$

\begin{cor}\label{cor4}
Let $A, B \in {\mathbb{C}}^{n\times n}$. If $s\geq max\{ind(AB), ind(BA)\}$, then $(AB)^s \sim (BA)^s $.
\end{cor}
\noindent
{\it Proof}. This is obvious by choosing ``$B = C$'' in Theorem~\ref{thm4}.      \hfill$\Box$   \\
\vskip-2mm
\indent
The following example illustrates Theorem~\ref{thm4} is a nontrivival generalization of~\cite[Corollary 2]{HARTWIG}.
\begin{exam}\rm
Let
$A=B=\begin{bmatrix}
         0 & 1 \\
         0 & 1 \\
       \end{bmatrix}$,
$C=\begin{bmatrix}
         0 & 1 \\
         1 & 0 \\
       \end{bmatrix} \in {\mathbb{C}^{2\times 2}}$.
Then $ABA=ACA$, while $B\neq C$. In this case, $ind(AC)=ind(BA)=1$. In view of Theorem~\ref{thm4}, $(AC)^s \sim (BA)^s$ for all $s\geq 1$.
Indeed,
$AC=\begin{bmatrix}
         1 & 0 \\
         1 & 0 \\
       \end{bmatrix}$ and
$BA=\begin{bmatrix}
         0 & 1 \\
         0 & 1 \\
       \end{bmatrix}$
are idempotent.
Let
$U=\begin{bmatrix}
         0 & 1 \\
         1 & 0 \\
       \end{bmatrix}$,
we have $(AC)^s =U^{-1}(BA)^s U $ for all $s\geq 1$.      \hfill$\Box$
\end{exam}

\end{document}